\documentclass[12pt,a4paper]{article}
\usepackage{amsmath,amssymb,amsfonts}
\def\Bbb{\mathbb}

\title{\bf  The rational part of a periodic continued fraction}

\author{Kurt Girstmair}
\date{}

\makeatletter
\let\@@maketitle=\maketitle
\def\maketitle{\def\thispagestyle##1{\relax}\@@maketitle}
\makeatother
\newtheorem{theorem}{Theorem}

\def\BE{\begin{equation}}
\def\EE{\end{equation}}
\def\BD{\begin{displaymath}}
\def\ED{\end{displaymath}}
\def\BA{\begin{array}}
\def\EA{\end{array}}
\def\BEA{\begin{eqnarray*}}
\def\EEA{\end{eqnarray*}}
\def\BI{\bibitem}

\def\Z{\Bbb Z}
\def\Q{\Bbb Q}

\def\phi{\varphi}
\def\EPS{\varepsilon}

\def\MB{\mbox}
\def\LD{\ldots}
\def\OV{\overline}

\def\BQ{``}
\def\EQ{'' }
\def\EQP{''}

\def\MN{\medskip\noindent}
\def\STOP{\hfill$\Box$}

\def\XQ{x_{\Q}}

\begin{document}
\maketitle

\begin{abstract}

\noindent

Let $x$ be a periodic continued fraction with the initial block $0$ and the repeating block $c_1,\ldots,c_n$. So $x$ is a quadratic irrational of
the form $x=a+\sqrt b$, where $a$, $b$ are rational numbers, $b>0$, $b$ not a square. The numbers $a$ and $\sqrt b$ are
uniquely determined by $x$. In general it is difficult to say what the influence of a certain digit of the repeating block on the appearance of $x$ is. We highlight a noteworthy exception from this rule.
Indeed, the magnitude of $2a$ is
essentially determined by the last digit $c_n$ of the repeating block, the fractional part of $2a$, however, is independent of $c_n$.  Of particular interest is the case $2a\in\Bbb Z$, which
occurs if, and only if, the sequence $c_1,\ldots,c_{n-1}$ is symmetric.

\end{abstract}

Continued fractions $x=[c_0,c_1,c_2,\LD]$ have aroused the interest of mathematicians for at least three centuries. They have many beautiful properties. But in general
they are, in some sense, hardly predictable. With the exception of the digit $c_0$ (which determines the magnitude of $x$), it is not easy to say what the influence of a certain digit $c_k$ on the
appearance of $x$ is. In this short note we present another noteworthy exception from this \BQ rule\EQP.

Let $x$ be an irrational real number. Then $x$ has an infinite continued fraction expansion $c_0,c_1, c_2 \LD$
The {\em convergents} $p_k/q_k$ of $x$ are defined in the usual way, namely,
\begin{eqnarray}
\label{1.4}
  p_{-1}=1,& p_{0}=c_0, \enspace p_k=c_kp_{k-1}+p_{k-2},&\nonumber\\
  q_{-1}=0,& q_{0}=1,   \enspace q_k=c_kq_{k-1}+q_{k-2}&
\end{eqnarray}
for $k\ge 1$. One also writes $p_k/q_k=[c_0,c_1,\LD,c_k]$ for $k\ge 0$. Note that $p_k$ and $q_k$ are relatively prime (see \cite[p. 3]{RoSz}).
The sequence $p_k/q_k$, $k\ge 0$, converges to $x$ (see \cite[p. 3]{RoSz}).
Accordingly, we write $x=[c_0,c_1,c_2,\LD]$.

We also need the $(n+1)$th {\em complete} quotient of $x$, which can be defined by
\BD
x_{n+1}=[c_{n+1},c_{n+2},\LD].
\ED
It satisfies
\BE
\label{1.6}
  x=\frac{x_{n+1}p_n+p_{n-1}}{x_{n+1}q_n+q_{n-1}}.
\EE
(see \cite[p. 4]{RoSz}).

A quadratic irrational is a real number $x$ of the form $x=a+\sqrt b$, where $a$, $b$ are rational numbers, $b>0$, $b$ not a square.
The root $\sqrt b$ may be positive or negative.
It is easy to see that $a$ and $\sqrt b$ are uniquely determined by $x$. We call $a$ the {\em rational part} of $x$ and
write $a=\XQ$. Similarly, $\sqrt b$ is called the {\em irrational part} of $x$.

Quadratic irrationals $x$ have infinite {\em periodic} continued fraction expansions, i.e.,
\BD
  x=[b_0,\LD,b_m,c_1,\LD,c_n,c_1,\LD,c_n,\LD],\enspace n\ge 1,
\ED
for which we also write $x=[b_0,\LD,b_m,\OV{c_1,\LD,c_n}]$ (see \cite[p. 41]{RoSz}). Here $b_0,\LD,b_m$ is called the {\em initial block}, $c_1,\LD,c_n$ the {\em repeating block} of $x$.
If there is no initial block, $x$ has the form $x=[\OV{c_1,\LD,c_n}]$ and is called {\em purely periodic} (see \cite[p. 39]{RoSz}).

Conversely, each periodic continued fraction expansion belongs to a quadratic irrational (see \cite[p. 40]{RoSz}).
In this note we consider quadratic irrationals that are
almost purely periodic, namely, of the form
\BE
\label{2.2}
  x=[0,\OV{c_1,\LD,c_n}].
\EE
Note that $\lfloor x\rfloor=0$, but $x\ne 0$, so $0<x<1$.
Since
$x=[0,[\OV{c_1,\LD,c_n}]]=0+1/[\OV{c_1,\LD,c_n}]$,
we see that $y=1/x=[\OV{c_1,\LD,c_n}]$ is purely periodic.

Our main objective is to describe the influence of the digit $c_n$ on the rational part $\XQ$ of $x$.
For this purpose we also consider the {\em fractional part} $\{a\}$ of a rational number $a$, which is defined by
\BD
  a-\{a\}\in\Z, 0\le \{a\}<1.
\ED

If $x$ is is as in (\ref{2.2}), the sequences $p_k$, $k\ge 0$, and $q_k$, $k\ge-1$, are increasing. This follows from (\ref{1.4}) by induction. Induction also gives $0\le p_k\le q_k$, $k\ge 0$.

\begin{theorem}
\label{t1}
Let $x$ be as in {\rm (\ref{2.2})}. Then
\BD
  2\XQ=-c_n+\EPS,
\ED
where $|\EPS|<1$ and $\EPS$  depends only of $c_1,\LD,c_{n-1}$. In particular, $\{2\XQ\}$ depends only of $c_1,\LD,c_{n-1}$.
\end{theorem}

\MN
{\em Remark.} The theorem says that the connection between $c_n$ and $2\XQ$ is remarkably simple. Whereas $-c_n$ basically determines the magnitude of $2\XQ$,
the fractional part of $2\XQ$ is independent of $c_n$. The influence of $c_n$ on the irrational part of $x$ is much more involved, as we will outline below.

\MN
{\em Proof of Theorem \ref{t1}.} As above, put $y=1/x=[\OV{c_1,\LD,c_n}]$.
Then we have
\BD
x=[0,c_1,\LD,c_n, y]=\frac{p_ny+p_{n-1}}{q_ny+q_{n-1}},
\ED
by (\ref{1.6}).
If we replace $y$ by $1/x$, we obtain
\BD
 x=\frac{p_n+p_{n-1}x}{q_n+q_{n-1}x}.
\ED
This gives rise to the quadratic equation
\BE
\label{2.3}
 q_{n-1}x^2+(q_n-p_{n-1})x-p_n=0,
\EE
which has the solutions $x=a+\sqrt b$, $x'=a-\sqrt b$. However, $x+x'=2a=2\XQ$, and a well-known fact about the roots of a quadratic equation, combined with (\ref{2.3}), says
\BD
 x+x'=\frac{p_{n-1}-q_n}{q_{n-1}}.
\ED
Now we use $q_n=c_nq_{n-1}+q_{n-2}$ and obtain
\BD
 2\XQ=-c_n+\frac{p_{n-1}-q_{n-2}}{q_{n-1}}.
\ED
We put
\BE
\label{2.4}
  \EPS=\frac{p_{n-1}-q_{n-2}}{q_{n-1}}.
\EE
Then $2\XQ=-c_n+\EPS$, and $\EPS$ depends only of $c_1,\LD,c_{n-1}$, not of $c_n$.
As we noted above, the convergents of $x$ satisfy $0\le p_k/q_k\le 1$ for all $k\ge 0$. In particular, $0\le p_{n-1}/q_{n-1}\le 1$. On the other hand,
the sequence $q_{-1}, q_{0},q_1,\LD$ (with $q_{-1}=0$ and $q_0=1$) is increasing, and so $0\le q_{n-2}/q_{n-1}\le 1$. This shows $|\EPS|\le 1$. If $\EPS=1$, then
$q_{n-2}=0$. Accordingly, $n=1$ and $p_{n-1}=0$, which gives $\EPS=0$, a contradiction. If $\EPS=-1$,
then $p_{n-1}=0$ and $n=1$ (observe $p_1=1$). Therefore, $q_{n-2}=0$ and $\EPS=0$, another contradiction. Hence $|\EPS|<1$ in all cases.
\\ \MB{ } \STOP

\MN
{\em Remark.} It is not difficult to describe $\{2\XQ\}$ in terms of $\EPS$. Indeed, since $|\EPS|<1$, (\ref{2.4}) shows $\{\EPS\}=\EPS$ if $p_{n-1}\ge q_{n-2}$,
and $\{\EPS\}=\EPS+1$, if $p_{n-1}<q_{n-2}$. Both cases are possible.

\MN
{\em Example.} Let $x=[0,\OV{1,2,2,3}]$. We have $n=4$ and $p_2/q_2=2/3$, $p_3/q_3=5/7$. Since $c_4=3$, we obtain $2\XQ=-3+\EPS$ with $\EPS=(5-3)/7=2/7=\{2\XQ\}$, by (\ref{2.4}).
Accordingly, $2\XQ=-19/7$. In order to find $x$ itself, we have to solve equation (\ref{2.3}), which reads $7x^2+19x-17=0$. We obtain $x=(-19 +\sqrt{837})/14$, where the root must be positive
since, otherwise, $x<0$. If we consider $x=[0,\OV{1,2,2,5}]$ instead, we already know that $2\XQ=-5+2/7=-33/7$, whereas the irrational part of $x$ turns out to be $\sqrt{1845}/14$.

\MN
{\em Remark.} An inspection of (\ref{2.3}) shows that the irrational part of $x$ is the square-root of a quadratic polynomial in $c_n$ whose coefficients are rational functions in $p_{n-1}$, $p_{n-2}$,
$q_{n-1}$, and $q_{n-2}$.

\MN
Of particular interest is the case when $2x_q\in\Z$, i.e., $\EPS=0$. Here we use the formula
\BD
  p_kq_{k-1}-q_kp_{k-1}=(-1)^{k+1},\enspace k\ge 0
\ED
(see \cite[p. 2]{RoSz}).
In our situation, it implies
\BE
\label{2.5}
  p_{n-1}q_{n-2}\equiv (-1)^n \mod q_{n-1}.
\EE
We have $\EPS=0$ if, and only if $p_{n-1}=q_{n-2}$. This is the case if $n=1$. If $n\ge 2$, we have $1\le p_{n-1}, q_{n-2}\le q_{n-1}$. Accordingly, $\EPS=0$ is the same as saying $p_{n-1}\equiv q_{n-2}\mod q_{n-1}$.
But then (\ref{2.5}) shows that $\EPS=0$ is equivalent to
\BE
\label{2.6}
p_{n-1}^2\equiv(-1)^n\mod q_{n-1}.
\EE
A classical theorem (see \cite[p. 28]{Pe}) says that (\ref{2.6}) is equivalent to the symmetry of $c_1,\LD,c_{n-1}$, i.e., $c_k=c_{n-k}$ for $k=1,\LD,n-1$ (this includes the case $n=1$).
Hence we obtain the following result.

\begin{theorem}
\label{t2}
In the setting of Theorem \ref{t1}, the following assertions are equivalent:
\\ {\rm (a)} $2\XQ\in\Z$.
\\ {\rm (b)} $2\XQ=-c_n$.
\\ {\rm (c)} The sequence $c_1,\LD,c_{n-1}$ is symmetric.
\end{theorem}

\MN
{\em Example.} Consider $x=[0,\OV{2,3,1,3,2,1}]$. We have $n=6$, $p_4/q_4=15/34$, $p_5/q_5=34/77$. Accordingly, $\EPS=(34-34)/77=0$ and $2\XQ=-1$.
Indeed, it turns out that $x=-1/2+\sqrt{39/44}$, where the square root is positive. So $x+1/2=\sqrt{39/44}$. Because the continued fraction expansions of $x$ and $x+1$ differ only in the digit $c_0$,
one might think that the continued fraction expansions of $x$ and
$x+1/2$ are also similar. This, however, is not the case since $\sqrt{39/44}=[0,1, \OV{16,11,1,3,2,3,1,11,16,2}]$. So we have another example of the \BQ unpredictability\EQ of continued
fractions we mentioned above.

\MN
{\em Remark.} The reader may find out whether something similar is true in the purely periodic case.


\vspace{0.5cm}
\noindent
Kurt Girstmair            \\
Institut f\"ur Mathematik \\
Universit\"at Innsbruck   \\
Technikerstr. 13/7        \\
A-6020 Innsbruck, Austria \\
Kurt.Girstmair@uibk.ac.at

\end{document}